\documentclass[12pt,a4paper, twoside]{article}

\usepackage{amsmath,amssymb,amsthm,amsfonts,mathrsfs,amscd,environ}
\usepackage{dsfont}
\usepackage{latexsym,enumerate,color,geometry,extarrows}
\usepackage{xcolor}
\usepackage{verbatim,fancyhdr,enumitem}

 \NewEnviron{ews}{%
\begin{equation}\begin{split}
  \BODY
\end{split}\end{equation}
}

\NewEnviron{ews*}{%
\begin{equation*}\begin{split}
  \BODY
\end{split}\end{equation*}
}

 \newtheorem{thm}{Theorem}[section]
 \newtheorem{lem}[thm]{Lemma}
 
 \newtheorem{exa}[thm]{Example}
 
 \newtheorem{defn}{Definition}[section]
 \newtheorem{rem}{Remark}[section]
 
 \numberwithin{equation}{section}

\def\dif{{\mathord{{\rm d}}}}
\def\no{\nonumber}

\def\mR{{\mathbb R}}
\def\mE{{\mathbb E}}

\def\mN{{\mathbb N}}

\def\mE{{\mathbb E}}
\def\mF{{\mathbb F}}

\def\mN{{\mathbb N}}

\def\mR{{\mathbb R}}

\def\sF{{\mathscr F}}

\def\sB{{\mathscr B}}
\def\sC{{\mathscr C}}

\def\sF{{\mathscr F}}

\def\sL{{\mathscr L}}
\def\sM{{\mathscr M}}

\def\sP{{\mathscr P}}

\def\bd{\begin{defn}}
\def\ed{\end{defn}}
\def\bl{\begin{lem}}
\def\el{\end{lem}}
\def\bt{\begin{thm}}
\def\et{\end{thm}}
\def\br{\begin{rem}}
\def\er{\end{rem}}

\allowdisplaybreaks

\title{{\bf Stabilization of stochastic McKean-Vlasov equations with feedback control based on discrete-time state observation }
}

\author{
{\bf Hao Wu$^{a)}$, Junhao Hu$^{a)}$,  Shuaibin Gao$^{b)},$ Chenggui Yuan$^{c)}$}\\
\footnotesize{$^{a)}$School of Mathematics and Statistics,
South-Central University For Nationalities}\\
\footnotesize{ Wuhan, Hubei 430000, P.R.China}\\
\footnotesize{Email: wuhaomoonsky@163.com},
\footnotesize{ junhaohu74@163.com}\\
\footnotesize{$^{b)}$ Department of Mathematics, Shanghai Normal University, Shanghai, 200234,  P.R.China}\\
\footnotesize{Email: shuaibingao@163.com}\\
\footnotesize{$^{c)}$Mathematic department, Swansea University, Bay campus, SA1 8EN, UK}\\
\footnotesize{Email: C.Yuan@Swansea.ac.uk}\\
}

\begin{document}

\maketitle

\begin{abstract}
In this paper, we  study the stability  of solutions of stochastic McKean-Vlasov equations (SMVEs) via feedback control based on discrete-time state observation. By using  a specific Lyapunov function,  the $H_{\infty}$ stability, asymptotic stability and  exponential stability in mean square for the solution of the controlled systems are obtained.
Since the distribution of  solution is  difficult to be observed, we study  the corresponding particle system which can be  observed for the feedback control.  We  prove that the exponential stability of control system  is equivalent to   the  the exponential stability of  the corresponding particle system. Finally,  an example is provided to show the effectiveness of the theory.
\end{abstract}\noindent

AMS Subject Classification (2020): \quad 60H10;\quad 93D15; \quad 60K35
\noindent

Keywords: McKean-Vlasov; Feedback control;  Stabilization; Discrete-time state observation; Distribution observation.

\section{Introduction}

 Stochastic differential equations (SDEs)   are widely used to model stochastic systems  in  different branches of science and industry.  The form of SDEs reads as follows:
 \begin{equation*}
\dif x(t)= f( x(t), t)\dif t
  +g(  x(t),t)\dif  B(t),  \, t\geq 0.
  \end{equation*}
     One of the popular  applications for SDEs   is   the feedback control of stochastic systems.  We refer the readers to \cite{ CZS, KKMM, RX1,  YMMH} and references therein. Since, some real SDEs are often unstable,  an interesting problem of the automatic control field is that for some given unstable SDEs, how can one  design an effective control function for the system to make the corresponding system be stable?  Among them, the  feedback control based on a continuous-time state observation is an efficient one, which has been used in establishing the mean-square exponential stabilization for a class of SDEs, see e.g. \cite{BHB, CFAF, FP1} and references therein. Since the method of continuous-time state observation is usually too expensive and not realistic in real lives, \cite{M1} proposed a more effective    state feedback control  which is  based on   discrete-time state observation and    is now widely studied.  It is obvious that  the state feedback control  based on   continuous-time observation requires one to observe the system all the time, while the state feedback control  based on   discrete-time state observation only requires one to observe the system in some discrete time. There are many results on this problem in the previous literatures (e.g. \cite{LMMY, SZQM}).  In particular, for an unstable stochastic system, it is very meaningful and important to design a feedback control with the form $u(\lfloor\frac{t}{\delta}\rfloor \delta)$ embedded into the drift part, where $\delta$ is the discrete-time observation gap.

On the other hand,  recently, many researchers are interested in  studying  the following equations, which are called stochastic Mckean-Vlasov equations (SMVEs):
\begin{equation*}
\begin{cases}
  &x(t)= x_{0}+\int^{t}_{0}f( x(s), \mu_{s})\dif s +g(   x(s), \mu_{s})\dif  B(s), t\in [t_{0}, \infty),\\
 &\mu_{t}=\sL(x(t)):= \mbox{the probability distribution of } x(t).
 \end{cases}
 \end{equation*}
Obviously, the coefficients involved depend not only on  the state process but also on its distribution. With contrast to the classical SDEs, SMVEs enjoy some essential features. The work on SMVEs was initiated by McKean \cite {M2}, who was inspired by Kac's Programme in Kinetic Theory \cite{K}. Sznitman \cite{S} investigated the existence and uniqueness of the results under a global Lipschitz condition.  Wang \cite{W} studied the existence of invariant probability measures for SMVEs. Govindan and Ahmed \cite{GA} studied the exponential stability of the solutions for a semilinear SMVEs under the Lipschitz condition and linear growth condition.   Ding and Qiao \cite{DQ1, DQ2} derived the  existence and uniqueness of the solution with non-Lipschitz condition and analyzed the stability of the solutions for SMVEs, respectively. Furthermore, in addition to the theoretical values, this kind of equations also has a lot of applied values in social science, economics, engineering, etc. (see  e.g. \cite{CD}).

To the best our knowledge,   there is  little study  on the stabilization of SMVEs with feedback control based on discrete-time state observation. It is clear that the controlled  Mckean-Vlasov system  includes discrete-time state observations as well as   its distribution observations while feedback control systems independent of distribution only need to observe state of the systems.  In this paper, we shall study  the stabilization problem by using the  feedback control with a discrete-time version: for an unstable McKean-Vlasov system, we aim to make the Mckean-Vlasov system stable by designing  a discrete-state and its distribution feedback control on this system. Our main contributions are as follows:
\begin{itemize}
\item[$\bullet$]We are the first to study feedback control problem for SMVEs  based on discrete-time state observation.

\item[$\bullet$] The Lyapunov functions used in this article not only contain state of the solution but also the distribution of the solution, while the previous  Lyapunov functions used in the state feedback control system only contain state of the solution. This is an  essential feature.

\item[$\bullet$] We  study the asymptotic stability and exponential stability in mean square of the solution for  SMVEs  based on discrete-time state observation.

\item[$\bullet$] The distribution of  analytical solution $x(t)$ is  difficult to be observed while the empirical distribution can be observed more easily. Thus, we further study  the corresponding particle system.  We show that the exponential stability of control system  is equivalent to   the  the exponential stability of the  corresponding particle system.
\end{itemize}

We close this part by giving our organization  in this article. In Section 2, we introduce some necessary notations, research objects and necessary assumptions. In Section 3, we aim to  study the stability of solutions to SMVEs via feedback control based on discrete-time state observation.  Then, an  example is presented to illustrate the theories.

 \section{Preliminaries}
\subsection{Notations}Throughout this paper,  let  $(\Omega, \sF, \mF, P)$  be  a complete probability space with filtration $\mF= \{\sF_{t}\}_{t\geq 0}$ satisfying the usual conditions(i.e., it is increasing and right continuous,  $\sF_{0}$ contains all $P$-null sets) taking along
 a standard  $m$-Brownian motion process  $B(t).$   If $x, y \in \mR^{d},$ we use $|x |$  to denote the Euclidean norm of
$x,$ and use  $\langle x, y\rangle$ or $xy$ to denote the Euclidean inner product. If  $A$ is a matrix, $A^{T}$  is the transpose  of $A,$  and   $|A |$ represents  $\sqrt{\mathrm{Tr} (AA^{T}).}$   Moreover, let $\lfloor a \rfloor$ be the integer parts of $a.$  For $\delta >0,$  set $\sigma_{t}=\lfloor\frac{t}{\delta}\rfloor \delta,$  where $\delta$ is the discrete-time observation gap.     Let $\sB(\mR^{d})$ be  the  Borel $\sigma-$algebra on $\mR^{d}$, $C(\mR^{d})$ denotes all continuous functions on $ \mR^{d}$ and
$C^{k}(\mR^{d})$ denotes all continuous functions on $ \mR^{d}$ with  continuous partial derivations
         of order up to $k$.
Let $\sP(\mR^{d})$ be the space  of all probability measures, and  $\sP_{p}(\mR^{d})$ denotes the space  of all probability measures defined on $\sB(\mR^{d}) $ with finite $p$th moment:
$$W_{p}(\mu):=\bigg(\bigg.\int_{\mR^{d}}|x|^{p}\mu(\dif x)\bigg)\bigg.^{\frac{1}{p}}<\infty.$$
For $ \mu, \nu\in \sP_{p}(\mR^{d})$, we define the
Wasserstein distance for $p\geq 1$  as follows: $$W_{p}(\mu, \nu):=\inf_{\pi \in \Pi(\mu, \nu)}\bigg\{\bigg. \int_{\mR^{d}\times \mR^{d}}|x-y|^{p}\pi(\dif x, \dif y)   \bigg\}\bigg.^{\frac{1}{p}}, $$
where $\Pi(\mu, \nu)$ is the family of all coupling for $\mu, \nu.$

Set
$\sM^{s}_{\lambda^{2}}(\mR^{d}):=\{m: m\, \mbox{is a signed measure on  } \mR^{d}\, \mbox{satisfying } $
{$\|\mu\|^{2}_{\lambda^{2}}=\int_{\mR^{d}}(1+|x|)^{2}|m|(\dif x)<\infty,$ \mbox{where}\, $|m|$ \,\mbox{is the total variation measure  of}\, $m$\},} and
$\sM_{\lambda^{2}}(\mR^{d})=\sM^{s}_{\lambda^{2}}(\mR^{d})\bigcap \sP(\mR^{d}).$
We put on $\sM_{\lambda^{2}}(\mR^{d})$ a topology induced by the  Wasserstein distance
$W_{2}(\cdot, \cdot).$

\subsection{  Lions Derivatives}
In this subsection,  we will give the definition of Loins derivative for $b: \sM_{\lambda^{2}}(\mR^{d}) \rightarrow  \mR$ with respect to a probability measure as introduced in  \cite{DQ2}.
\bd
We say that  $U: \sM_{\lambda^{2}}(\mR^{d}) \rightarrow  \mR$ is differential at $\mu \in \sM_{\lambda^{2}}(\mR^{d}),$ if there exists some $X\in L^{2}(\Omega, \sF,P; \mR^{d})$ such that $\mu=\sL(X)$ and the function $\tilde{U}:L^{2}(\Omega, \sF,P; \mR^{d}) \rightarrow \mR $ given by $\tilde{U}(X):=U(\sL(X))$ is Fr\'{e}chet differentiable at $X.$
\ed
We recall that $\tilde{U}$ is Fr\'{e}chet differentiable at $X$ means that there exists a continuous mapping $D\tilde{U}(X): L^{2}(\Omega, \sF,P; \mR^{d}) \rightarrow \mR$ such that for any $Y\in L^{2}(\Omega, \sF,P; \mR^{d})$
\begin{align*}
\tilde{U}(X+Y)-\tilde{U}(X)=D\tilde{U}(X)(Y)+o(|Y|_{L^{2}}), \mbox{as}\,\ |Y|_{L^{2}}\rightarrow 0.
\end{align*}
Due to $D\tilde{U}(X)\in L(L^{2}(\Omega, \sF,P; \mR^{d});\mR),$ by Riesz representation theorem, there exists  a $P$-a.s. unique variable  $Z\in L^{2}(\Omega, \sF,P; \mR^{d})$ such that for any $Y\in L^{2}(\Omega, \sF,P; \mR^{d})$
$$ D\tilde{U}(X)(Y)=\langle Y, Z\rangle_{L^{2}}= \mE[Y Z]. $$
Cardaliaguet \cite{C} showed  that there exists a Borel measurable function  $h:\mR^{d}\rightarrow \mR^{d}$ which only depends on the distribution $\sL(X)$ rather that $X$ itself such that $Z=h(X).$ Thus, for $X\in L^{2}(\Omega, \sF,P; \mR^{d}),$
$$U(\sL(Y))-U(\sL(X))= \mE[h(X)(Y-X)]+o(|Y-X|_{L^{2}}).$$
We call $\partial_{\mu} U(\sL(X))(y):=h(y), y\in \mR^{d}$ the L-derivative of $U$ at $\sL(X), X\in L^{2}(\Omega, \sF,P; \mR^{d}).$

 Let  $C^{1}(\sM_{\lambda^{2}}(\mR^{d}))$ denote all functions $U: \sM_{\lambda^{2}}(\mR^{d}) \rightarrow  \mR$ such that $\partial_{\mu} U: \sM_{\lambda^{2}}(\mR^{d})\times \mR^{d}\rightarrow \mR^{d}$ is continuous. Let $C^{1;1}_{b}(\sM_{\lambda^{2}}(\mR^{d}))$ be all functions $U\in C^{1}(\sM_{\lambda^{2}}(\mR^{d}))$ such that $\partial_{\mu} U $ is bounded and Lipschitz continuous, i.e., there exists a positive constant $C$ such that
 \begin{itemize}
\item[(i)] $|\partial_{\mu} U(\mu)(x)|\leq C$ for any $\mu\in  \sM_{\lambda^{2}}(\mR^{d}), x\in \mR^{d}.$
\item[(ii)] $|\partial_{\mu} U(\mu)(x)- \partial_{\mu} U(\nu)(y)|^{2}\leq C(W^{2}_{2}(\mu, \nu)+|x-y|^{2}),\mu, \nu \in \sM_{\lambda^{2}}(\mR^{d}), x, y \in \mR^{d}.  $
\end{itemize}

We need more definitions:
\begin{description}
\item[(1)] The function $U $ is said to be in $C^{2}(\sM_{\lambda^{2}}(\mR^{d}))$  if for any $\mu\in \sM_{\lambda^{2}}(\mR^{d}),$ $U\in C^{1}(\sM_{\lambda^{2}}(\mR^{d})), $ $\partial_{\mu}U(\mu)(\cdot)$ is differentiable and its derivative $\partial_{y}\partial_{\mu}U: \sM_{\lambda^{2}}(\mR^{d})\times \mR^{d}\rightarrow \mR^{d}\otimes \mR^{d}$ is continuous.

\item[(2)]
The function $U $ is said to be in $C^{2;1}_{b}(\sM_{\lambda^{2}}(\mR^{d}))$  if $U\in C^{2}(\sM_{\lambda^{2}}(\mR^{d}))\bigcap  C^{1;1}_{b}$ $(\sM_{\lambda^{2}}(\mR^{d}))$ and its derivative $\partial_{y}\partial_{\mu}U$ is bounded and Lipschitz continuous.

\item[(3)]
The function $\Psi $ is said to be in $C^{2,2}_{b}(\mR^{d}\times\sM_{\lambda^{2}}(\mR^{d}))$  if for any $x\in \mR^{d},$ $\Psi(x,\cdot)\in C^{2}(\sM_{\lambda^{2}}(\mR^{d}))$ and for any $\mu\in \sM_{\lambda^{2}}(\mR^{d})$, $\Psi(\cdot,\mu)\in C^{2}(\mR^{d}).$

\item[(4)]
The function $\Psi $ is said to be in $\sC(\mR^{d}\times\sM_{\lambda^{2}}(\mR^{d}))$  if $\Psi\in C^{2,2}_{b}(\mR^{d}\times\sM_{\lambda^{2}}(\mR^{d}))$ and for any compact set $\mathcal{K}\subset \mR^{d}\times \sM_{\lambda^{2}}(\mR^{d}),$
 $$\sup_{(x,\mu)\in \mathcal{K}}\int_{\mR^{d}}(|\partial_{y}\partial_{\mu}\Psi (x,\mu)(y)|^{2}+|\partial_{\mu}\Psi (x,\mu)(y)|^{2})\mu(\dif y) <\infty.$$

\item[(5)]
 $\Psi \in C^{2,2;1}_{b}(\mR^{d}\times \sM_{\lambda^{2}}(\mR^{d}))$ means that
 \begin{itemize}
\item[(i)] $\Psi$ is bicontinuous  in $(x,\mu).$
\item[(ii)] For any $x\in \mR^{d},$ $\Psi(x,\cdot)\in C^{2;1}_{b}(\sM_{\lambda^{2}}(\mR^{d}))$ and for any $\mu\in \sM_{\lambda^{2}}(\mR^{d})$, $\Psi(\cdot,\mu)\in C^{2}(\mR^{d}).$
\item[(iii)] For any $\mu\in \sM_{\lambda^{2}}(\mR^{d}),$ $\partial_{x_{i}}\Psi (\cdot,\mu)$ is bounded.
\end{itemize}

\item[(6)]  Let $\sC_{+}(\mR^{d}\times\sM_{\lambda^{2}}(\mR^{d}))$ be  a set of all functions $\Psi \in \sC(\mR^{d}\times\sM_{\lambda^{2}}(\mR^{d}))$ such that $\Psi >0.$
 $ C^{2,2;1}_{b,+}(\mR^{d} \times \sM_{\lambda^{2}}(\mR^{d}))$ denotes all functions $\Psi \in C^{2,2;1}_{b}(\sM_{\lambda^{2}}(\mR^{d}))$ with $\Psi \geq0.$

\end{description}

\subsection{The It\^o formula}

Consider the following equations:
\begin{equation}\label{1}
\begin{cases}
 &\dif y(t)= (f( y(t), \rho_{t})+u( y(t), \rho_{t}))\dif t +g(  y(t), \rho_{t})\dif  B(t), t\in [0, \infty),\\
 & y(0)=x_{0},
 \end{cases}
 \end{equation}
and
\begin{equation}\label{2}
\begin{cases}
 &\dif x(t)= (f( x(t), \mu_{t})+u( x(\sigma_{t}), \mu_{\sigma_{t}}))\dif t +g(  x(t), \mu_{t})\dif  B(t), t\in [0, \infty),\\
 & x(0)=x_{0},
 \end{cases}
 \end{equation}
where  $x_{0}\in \mR^{d},$  and $\rho_t$ and $\mu_t$ are the distributions of $y(t)$ and $x(t)$, respectively.  Moreover,  $f, u: \mR^{d}\times \sM_{\lambda^{2}}(\mR^{d})\rightarrow \mR^{d}, g:  \mR^{ d}\times \sM_{\lambda^{2}}(\mR^{d})\rightarrow \mR^{d\times m}.$
 In  Eq.\eqref{2}, one can see that the control function $u(x(\sigma_{t}), \mu_{\sigma_{t}})$  only depends on  the state at discrete times $0, \delta, 2\delta, \cdots.$ Moreover, we assume that $f({\bf0}, \delta_{\bf0})=0, g({\bf0}, \delta_{\bf0})=0, \mu({\bf0}, \delta_{\bf0})=0.$
For the existence and uniqueness of Eq.\eqref{1} and Eq.\eqref{2}, we assume that:
\begin{itemize}
\item[(H1)] Suppose that $f, g, u$  satisfy the following Lipschitz condition,  i.e., there exist  positive constants $L_{i}, i=1,2,3$ such that
$$|f(x,\mu)-f(y,\nu)|^{2}\leq L_{1}(|x-y|^{2}+W^{2}_{2}(\mu, \nu)), $$
$$|g({x},\mu)-g({y},\nu)|^{2}\leq L_{2}(|x-y|^{2}+W^{2}_{2}(\mu, \nu)), $$
$$|u(x,\mu)-u(y,\nu)|^{2}\leq L_{3}(|x-y|^{2}+W^{2}_{2}(\mu, \nu)), $$
for all $x,y\in \mR^{d}, \mu, \nu \in \sM_{\lambda^{2}}(\mR^{d}).$
\end{itemize}

By Theorem 3.1 in  \cite{HYY},  under the condition $(\mathrm{H1}),$  Eq.\eqref{1} and Eq.\eqref{2} have unique solutions, respectively.
We now introduce the following operators.
\bd
For $V\in  \sC(\mR^{d}\times\sM_{\lambda^{2}}(\mR^{d})),$  the operator $LV: \mR^{d}\times \sM_{\lambda^{2}}(\mR^{d})\rightarrow \mR$ for Eq.\eqref{1} is defined
 \begin{align}\label{II}
&LV(x, \rho)=\partial_{x}V(x, \rho)f(x, \rho)
 + \partial_{x}V(x, \rho)u(x, \rho)\no\\
 &+ \frac{1}{2}{\rm tr}( g^{T}(x, \rho)V_{xx}(x, \rho)g(x, \rho) )
 + \int_{\mR^{d}}\partial_{\rho}V(x, \rho)(y)f(y,\rho)\rho(\dif y)\no\\
 &+ \int_{\mR^{d}}\partial_{\rho}V(x, \rho)(y)u(y,\rho)\rho(\dif y)+\frac{1}{2}\int_{\mR^{d}}{\rm tr}(g^{T}(y, \rho)\partial_{y} \partial_{\rho}V(x, \rho)(y)g(y, \rho)) \rho(\dif y).
\end{align}
\ed

\bd
Let $\xi$ and $\eta$ be  two random variables  whose distributions  are $\mu$ and $\nu$,  respectively.  Let  the joint distribution of  $(\xi, \eta)$ be $F_{\xi, \eta}(z, \bar z).$  For $V\in  \sC(\mR^{d}\times\sM_{\lambda^{2}}(\mR^{d})),$  the  operator $\mathbb{L}V: \mR^{d}\times \sM_{\lambda^{2}}(\mR^{d})\times\mR^{d}\times \sM_{\lambda^{2}}(\mR^{d})\rightarrow \mR $ for Eq.\eqref{2} is defined by
 \begin{align*}
\mathbb{L}V(x, \mu,y, \nu)&=\partial_{x}V(x, \mu)f(x, \mu)
 + \partial_{x}V(x, \mu)u(y, \nu)\no\\
 &+ \frac{1}{2}{\rm tr}( g^{T}(x, \mu)V_{xx}(x, \mu)g(x, \mu) )
 + \int_{\mR^{d}}\partial_{\mu}V(x, \mu)(z)f(z,\mu)\mu(\dif z)\no\\
 &+ \int_{\mR^{d}}\int_{\mR^{d}}\partial_{\mu}V(x, \mu)(z)u(\bar{z},\nu)F_{\xi, \eta}(\dif z, \dif \bar z)\no\\
 &+\frac{1}{2}\int_{\mR^{d}}{\rm tr}(g^T(z, \mu)\partial_{z} \partial_{\mu}V(x, \mu)(z) g(z, \mu) ) \mu(\dif z).
\end{align*}
\ed
The It\^o's  formula has been established in \cite{DQ2,HSS} for equation \eqref{1}, we cite it as the following lemma.
\bl
Assume $(\mathrm{H1})$ and  $V\in C^{1, 2, (1, 1)}(\mR^{d}\times\sM_{\lambda^{2}}(\mR^{d}))$. Then it holds that
\begin{equation*}\label{Ito1}
V(y(t), \rho_t)-V(y(0), \rho_0)=\int_0^t L V(y(s), \rho_s)\dif s+ \int_0^t V_x(y(s), \rho_s)g(y(s), \rho_s)\dif B(s).
\end{equation*}
\el
Since the feedback control in  Eq.\eqref{2}  depends on the discrete time, we need to develop an It\^o's formula for this equation. One can find the It\^o formula for Eq. \eqref{1} in \cite{DQ2,HSS}.

We need more notations to formulate the It\^o formula.
Assume that $(\tilde{\Omega}, \tilde{\sF}, \tilde{\mF}=\{\tilde{\sF}_{t}\},   \tilde{P})$ is another probability space taking along a $m-$dimensional Brownian motion $\{\tilde{B}(t)\}_{t\geq 0}$. Consider the following equation:
\begin{equation}\label{+2}
\begin{cases}
 &\dif \tilde{x}(t)= (f( \tilde{x}(t), \tilde{\mu}_{t})+u( \tilde{x}(\sigma_{t}), \tilde{\mu}_{\sigma_{t}}))\dif t +\tilde{g}(  \tilde{x}(t), \tilde{\mu}_{t})\dif  \tilde{B}(t), t\in [0, \infty),\\
 & \tilde{x}(0)=x_{0},
 \end{cases}
 \end{equation}
where  $x_{0}\in \mR^{d},$ and $\tilde{\mu}_{t}$ denotes the distribution of $\tilde{x}(t).$
By the weak uniqueness,  it holds that $\{x(t)\}_{t\geq 0}$ and $\{\tilde{x}(t)\}_{t\geq 0}$ are identical in probability law. Furthermore,  denote by $\tilde{\mE}[\cdot]$ the expectation under $\tilde{P}$.

 We now present  the It\^o formula for \eqref{2}, which is   an  extension of proposition 2.9 in \cite{DQ2}.
\bl
Let $V\in \sC(\mR^{d}\times\sM_{\lambda^{2}}(\mR^{d}))$  and  the assumption  $(\mathrm{H1})$ hold. Then one has that
\begin{align}\label{a1}
 V(x(t), \mu_{t})&= V(x(0), \mu_{0})
 + \int^{t}_{0}\mathbb{L}V(x(s), \mu_{s},x(\sigma_{s}), \mu_{\sigma_{s}})\dif s\no\\
 &+\int^{t}_{0} \partial_{x}V(x(s), \mu_{s})g(x(s), \mu_{s})\dif B(s).
\end{align}
\el

\begin{proof} Let $x(t)$ be the unique solution of \eqref{2}.
By  H\"{o}lder's inequality and BDG's   inequality,  it holds that for any $T>0$ and $0\leq t \leq T,$
\begin{align*}
 \mE\sup_{0\leq s\leq t}|x(s)|^{2}&\leq 3\mE|x_0|^{2}+3\mE\sup_{0\leq s\leq t}\bigg|\bigg.\int^{s}_{0} f( x(r),\mu_{r})+u(  x(\sigma_{r}), \mu_{\sigma_{r}})\dif r \bigg|\bigg.^{2}\\
 &+3\mE\sup_{0\leq s\leq t}\bigg|\bigg.\int^{s}_{0} g(  {x}(r), \mu_{r})\dif B(r) \bigg|\bigg.^{2}\\
 &\leq 3\mE|x_0|^{2}+6T\mE\int^{t}_{0} |f({x}(s), \mu_{s})|^{2}\dif s+6T\mE\int^{t}_{0} |u( x(\sigma_{s}), \mu_{\sigma_{s}})|^{2}\dif s \\
 &  \quad\quad\quad+12\mE\int^{t}_{0} |g(  {x}(s), \mu_{s})|^{2}\dif s\\
 &\leq 3\mE|x_0|^{2}+\mE\int^{t}_{0}[(6TL_{1}+12L_{2})|{x}(s)|^{2}+6TL_{3}| x(\sigma_{s})|^{2}\\
 &\quad\quad\quad+(6TL_{1}+12L_{2})W^{2}_{2}(\mu_{s}, \delta_{0})+6TL_{3}W^{2}_{2}(\mu_{\sigma_{s}}, \delta_{0})]\dif s.
\end{align*}
This, together with $W_{2}^{2}(\mu_{s}, \delta_{0})\leq \mE|{x}(s)|^{2}$,    yields that
\begin{align*}
 \mE\sup_{0\leq s\leq t}|{x}(s)|^{2}&\leq
 3\mE|x_0|^{2}+ (12TL_{1}+24L_{2}+12TL_{3})\mE\int^{t}_{0}\sup_{0\leq u\leq s}|{x}(u)|^{2} \dif s.
\end{align*}
From Gronwall's  formula, we get
 \begin{align}\label{a2}
 \mE\sup_{0\leq s\leq T}|{x}(s)|^{2}\leq  3\mE|\xi(0)|^{2}e^{6(2TL_{1}+2TL_{3}+4L_{2})T}.
\end{align}
Using similar method, one can derive that
  \begin{align}\label{a3}
 \mE\sup_{0\leq s\leq T}|  x(s)|^{4}\leq C,
\end{align}
where $C$ is a constant depending on $L_{1},  L_{2}, L_{3},T.$
By \eqref{a2},  \eqref{a3} and (H1), we have
 \begin{align}\label{ba}
 &\mE\int^{T}_{0}(|f( x(s),  \mu_{s})|^{2}+ |u(  x(\sigma_{s}), \mu_{\sigma_{s}})|^2+|g( x(s), \mu_{s})|^{4})\dif s\no \\
 &\leq (L_{1}+L_{3}) \mE\int^{T}_{0}(| x(s)|^{2}+| x(\sigma_{s})|^2+W_{2}^{2}( \mu_{s}, \delta_{0})+W_{2}^{2}( \mu_{\sigma_{s}}, \delta_{0}))\dif s\no\\
 &+ 2L^{2}_{2}\mE\int^{T}_{0}(| x(s)|^{4}+W^{4}( \mu_{s}, \delta_{0}))\dif s< \infty.
\end{align}
Let
$$
\tilde b_t=(f( \tilde x(t), \tilde \mu_{t})+u( \tilde x(\sigma_{t}), \tilde \mu_{\sigma_{t}})), \, \, \tilde \sigma_t=g(  \tilde x(t), \tilde \mu_{t}).
$$
The Eq.\eqref{+2} can be written by
\begin{align*}
d\tilde x(t)=\tilde b_t \dif  t+\tilde \sigma_t \dif B(t).
\end{align*}
Since $x(t)$ and $\tilde x(t)$ have the same distribution,  from \eqref{ba} one can see that
$$
\tilde \mE \int_0^t[|\tilde b_s|^2+|\tilde \sigma_s|^4]\dif t<\infty.
$$
Fix $x$ and   set $h(\mu)=V(x, \mu).$     It follows from \cite[Proposition A.6]{HSS}  that
\begin{align*}
h(\tilde \mu_{t})- h(\tilde  \mu_{0})&=\int^{t}_{0}{\mE}\bigg[\tilde b_s\partial_{\mu}V({x}, \tilde \mu_{s})(\tilde {x}(s))+\frac{1}{2}{\rm tr}[\tilde \sigma^*_s\partial_{x} \partial_{\mu}V( x, \tilde \mu_{s}) )(\tilde {x}(s)) \tilde \sigma_s]\bigg]\dif s,\no\\
&=\int^{t}_{0}{\mE}\bigg[\Big(f(\tilde {x}(s), \tilde \mu_{s})+u(  \tilde x(\sigma_{s}), \tilde \mu_{\sigma_{s}})\Big)\partial_{\mu}V({x}, \tilde \mu_{s})(\tilde {x}(s))\no\\
& + \frac{1}{2}{\rm tr}\Big(g^*(\tilde {x}(s), \tilde \mu_{s})\partial_{x} \partial_{\mu}V( x, \tilde \mu_{s}) )(\tilde {x}(s)) g(\tilde {x}(s), \tilde \mu_{s}) \Big)  \bigg]\dif s\no\\
&=:\int^{t}_{0}M({x}, \tilde \mu_{s})\dif s.
\end{align*}
Now, set $\overline{V}(x, t)=V( x, \tilde \mu_{t}).$ Thus, we have
$$\partial_{t}\overline{V}(x, t)=M({x},  \tilde \mu_{t}).$$
Applying It\^{o}'s formula  \cite[Proposition A.8]{HSS} to ${V}(x(t), \tilde\mu_t)$ and noting that $\mu_t=\tilde \mu_t$, we derive that
\begin{align*}
&{V}(x(t), \tilde \mu_t)-{V}(x(0), \tilde \mu_0)
=
\overline{V}(x(t), t)- \overline{V}(x(0), 0)\\
&=\int^{t}_{0}\bigg[V_x(x(s), \mu_{s})(x(s))\big[f(x(s), \mu_{s})+u(x(\sigma_s),  \mu_{\sigma_s})\big]\no\\
& + \frac{1}{2}{\rm trace}\big[g^T(x(s), \mu_{s})V_{xx}(x(s), \mu_{s})(z) g(x(s), \mu_{s})\big] +M(x(s), \tilde\mu_{s}) \bigg]\dif s\no\\
& + \int^{t}_{0}f(x(s), \mu_{s})\partial_{x}V(x(s), \mu_{s})g(x(s), \mu_{s})\dif B(s).
\end{align*}
 The desired assertion \eqref{a1} holds.
\end{proof}

\section{Asymptotic stability and exponential stability in mean square }

In order to study the asymptotic stability and the exponential  stability in mean square, we impose the following assumption:

\begin{itemize}
\item[(H2)] Assume  there exist $V\in \sC(\mR^{d}\times\sM_{\lambda^{2}}(\mR^{d}))$, and    four constants $\lambda_{1}>0, \lambda_{2}>0, \gamma_{1}>0, \gamma_{2}\geq 0$ such that
\begin{align*}\int_{\mR^{d}}&LV(x,\mu)\mu(\dif x)+\lambda_{1}\int_{\mR^{d}}|V_{x}(x,\mu)|^{2}\mu(\dif x)+ \lambda_{2}\int_{\mR^{d}}|\partial_{\mu}V(x,\mu)(y)|^{2}\mu(\dif y)\\
&\leq -\gamma_{1}\int_{\mR^{d}}V(x,\mu)\mu(\dif x)  +\gamma_{2}.
\end{align*}
\end{itemize}

The following two results are about the  asymptotic stability  of the solutions for Eq. \eqref {2}.
\bl
Let   $(\mathrm{H1})-(\mathrm{H2})$  hold and assume further that there exists a positive constant $c_1$ such that $c_{1}\int_{\mR^{d}}|x|^{2}\mu(\dif x)\leq \int_{\mR^{d}}V(x,\mu)\mu(\dif x).$ If $\delta>0$ is sufficiently small  such that
$$\gamma_{1}c_{1}-8L_{1}\theta\delta^{2}-8L_{3}\theta\delta^{2}-2L_{2}\theta\delta\geq 0 ,\delta \leq \sqrt{\frac{1}{8L_{3}}},$$
then the control system \eqref{2} is $H_{\infty}-$stable, i.e.,
\begin{align}\label{91Y}
\limsup_{t\rightarrow\infty}\frac{1}{t}\mE\int^{t}_{0}|x(s)|^{2}\dif s\le  \gamma_{2},
\end{align}
for all initial data $x_{0}\in \mR^{n}.$  Moreover, if $\gamma_{2}=0, $ we have
\begin{align}\label{91}
\mE\int^{\infty}_{0}|x(t)|^{2}\dif t< \infty.
\end{align}
\el
\begin{proof}
We divide the proof into  two parts.

(i)  {\it We construct the following Lyapunov  functional which depends on the segment process $x_{t}:=\{x(t+r); -\delta \leq r\leq 0\}$ with $x(r)=x_{0}\in \mR^{d}, -\delta\leq r \leq 0.$ } That is: Let
\begin{align}\label{7}
 \widetilde{V}(x_{t}, \mu_{t})= V(x(t), \mu_{t})+\theta\int^{t}_{t-\delta}&\int^{t}_{r}[\delta|f(x(s),\mu_{s})\no\\
 &+u(x(\sigma_{s}), \mu_{\sigma_{s}})|^{2}+|g(x(s), \mu_{s})|^{2}]\dif s\dif r, \, t\geq 0,
\end{align}
where $\theta$ is a positive constant to be determined later.  Applying It\^{o}'s formula to $\widetilde{V}(x_t, \mu_{t})$ and noting that
\begin{align*}
\dif & \bigg(\bigg. \theta\int^{t}_{t-\delta}\int^{t}_{r}[\delta|f(x(u), \mu_{u})
 +u(x(\sigma_{u}), \mu_{\sigma_{u}})|^{2}+|g(x(u), \mu_{u})|^{2}]\dif u\dif r\bigg)\bigg.\\
 & =\bigg[\bigg. \theta\delta[\delta|f(x(t), \mu_{t})+u(x(\sigma_{t}), \mu_{\sigma_{t}})|^{2}+|g(x(t), \mu_{t})|^{2}]\no\\
 &-\theta\int^{t}_{t-\delta}[\delta|f(x(r), \mu_{r})
 +u(x(\sigma_{r}), \mu_{\sigma_{r}})|^{2}+|g(x(r), \mu_{r})|^{2}]\dif r\bigg]\bigg. \dif t,
\end{align*}  we get
\begin{align}\label{8}
\dif  \widetilde{V}(x_{t}, \mu_{t})= \mathscr{L}{V}(x_{t}, \mu_{t})+ \dif M(t),
\end{align}
where $M(t)$ is a martingale and
 \begin{align}\label{9}
\mathscr{L}{V}(x_{t}, \mu_{t})&=\mathbb{L}V(x(t), \mu_t, x(\sigma_t), \mu_{\sigma_t})
 +\theta\delta[\delta|f(x(t), \mu_{t})+u(x(\sigma_{t}), \mu_{\sigma_{t}})|^{2}+|g(x(t), \mu_{t})|^{2}]\no\\
 &-\theta\int^{t}_{t-\delta}[\delta|f(x(r), \mu_{r})
 +u(x(\sigma_{r}), \mu_{\sigma_{r}})|^{2}+|g(x(r), \mu_{r})|^{2}]\dif r.
\end{align}

(ii) {\it We are going to prove \eqref{91Y} and \eqref{91}}.
From \eqref{II} and \eqref{9}, we get
\begin{align}\label{12}
 \mathscr{L}{V}(x_{t}, \mu_{t})&=LV(x(t), \mu_{t})-\partial_{x}V(x(t), \mu_{t})(u(x(t), \mu_{t})-u(x(\sigma_{t}), \mu_{\sigma_{t}})) \no\\
 &- \int_{\mR^{d}}\int_{\mR^{d}}\partial_{\mu}V(x(t), \mu_{t})(y)(u(y,\mu_{t})-u(\bar{y},\mu_{\sigma_t}))F_{x(t), x(\sigma(t)) }(\dif y, \dif \bar{y})\no\\
 &+\theta\delta[\delta|f(x(t), \mu_{t})+u(x(\sigma_{t}), \mu_{\sigma_{t}})|^{2}+|g(x(t), \mu_{t})|^{2}]\no\\
 &-\theta\int^{t}_{t-\delta}[\delta|f(x(r), \mu_{r})
 +u(x(\sigma_{r}), \mu_{\sigma_{r}})|^{2}+|g(x(r), \mu_{r})|^{2}]\dif r.
\end{align}
By Young's inequality, we have
\begin{align}\label{13}
-\partial_{x}&V(x(t), \mu_{t})(u(x(t), \mu_{t})-u(x(\sigma_{t}), \mu_{\sigma_{t}}))\no\\
&\leq \lambda_{1}|\partial_{x}V(x(t), \mu_{t})|^{2}+\frac{L_{3}}{4\lambda_{1}}|x(t)-x(\sigma_{t})|^{2} +\frac{L_{3}}{4\lambda_{1}}W_{2}^{2}(\mu_{t},\mu_{\sigma_{t}}),
\end{align}
and
\begin{align}\label{14}
- \int_{\mR^{d}}&\int_{\mR^{d}}\partial_{\mu}V(x(t), \mu_{t})(y)(u(y,\mu_{t})-u(\bar{y},\mu_{\sigma_t}))F_{x(t), x(\sigma(t)) }(\dif y, \dif \bar{y})\no\\
&\leq \lambda_{2} \int_{\mR^{d}}|\partial_{\mu}V(x(t), \mu_{t})(y)|^{2}\mu_{t}(\dif y)+\frac{L_{3}}{4\lambda_{2}}\mE[|x(t)-x(\sigma_{t})|^{2}]+\frac{L_{3}}{4\lambda_{2}}W_{2}^{2}(\mu_{t},\mu_{\sigma_{t}}).
\end{align}
Since $ W_{2}^{2}(\mu_{t},\delta_{0})\leq \mE|x(t)|^{2},$ it follows from \eqref{12}, \eqref{13}, \eqref{14}   and (H1) that
\begin{align}\label{15}
 \mathscr{L}&{V}(x_{t}, \mu_{t})\le LV(x(t), \mu_{t})+ \lambda_{1}|\partial_{x}V(x(t), \mu_{t})|^{2}+\lambda_{2}\int_{\mR^{d}}|\partial_{\mu}V(x(t), \mu_{t})(y)|^{2}\mu_{t}(\dif y)\no\\
 &+(4L_{1}\theta\delta^{2}+4L_{3}\theta\delta^{2}+L_{2}\theta\delta)|x(t)|^{2}+ (4L_{1}\theta\delta^{2}+4L_{3}\theta\delta^{2}+L_{2}\theta\delta)W_{2}^{2}(\mu_{t}, \delta_{0})\no\\
 &+(\frac{L_{3}}{4\lambda_{1}}+2L_{3}\theta\delta^{2}+\frac{L_{3}}{4\lambda_{2}})|x(t)-x(\sigma_{t})|^{2}+(\frac{L_{3}}{4\lambda_{1}}
 +2C_{3}\theta\delta^{2}+\frac{L_{3}}{4\lambda_{2}})W_{2}^{2}(\mu(t), \mu(\sigma_{t}))\no\\
 &-\theta\int^{t}_{t-\delta}[\delta|f(x(r), \mu_{r})
 +u(x(\sigma_{r}), \mu_{\sigma_{r}})|^{2}+|g(x(r), \mu_{r})|^{2}]\dif r.
\end{align}
Noting that
$W_{2}^{2}(\mu_{t}, \mu_{\sigma_{t}})\leq \mE|x(t)-x(\sigma_{t})|^{2}, $
 \eqref{15} and $(\mathrm{H2})$, we obtain
\begin{align}\label{16}
&\mE \mathscr{L}{V}(x_{t}, \mu_{t})=\int_{\mR^d}\mathscr{L} {V}(x, \mu_{t})\mu_{t}(\dif x)\le-\lambda_{4}\mE|x(t)|^{2}+\gamma_{2}\no\\
  &+(\frac{L_{3}}{4\lambda_{1}}+2L_{3}\theta\delta^{2}+\frac{L_{3}}{4\lambda_{2}})\mE|x(t)-x(\sigma_{t})|^{2}+(\frac{L_{3}}{4\lambda_{1}}
 +2L_{3}\theta\delta^{2}+\frac{L_{3}}{4\lambda_{2}})\mE|x(t)-x(\sigma_{t})|^{2}\no\\
 &-\theta\mE\int^{t}_{t-\tau}[\delta|f(x(r), \mu_{r})
 +u(x(\sigma_{r}), \mu_{\sigma_{r}})|^{2}+|g(x(r), \mu_{r})|^{2}]\dif r\no\\
 &\leq-\lambda_{4}\mE|x(t)|^{2}+\gamma_{2}\no\\
  &+2(\frac{L_{3}}{4\lambda_{1}}+2L_{3}\theta\delta^{2}+\frac{L_{3}}{4\lambda_{2}})\mE|x(t)-x(\sigma_{t})|^{2}\no\\
 &-\theta\mE\int^{t}_{t-\delta}[\delta|f(x(r), \mu_{r})
 +u(x(\sigma_{r}), \mu_{\sigma_{r}})|^{2}+|g(x(r), \mu_{r})|^{2}]\dif r.
\end{align}
where $\lambda_{4}=\gamma_{1}c_{1}-8L_{1}\theta\delta^{2}-8L_{3}\theta\delta^{2}-2L_{2}\theta\delta.$
Noting that $t-\sigma_{t}\leq \delta,$  we have
\begin{align*}
 \mE|x(t)-x(\sigma_{t})|^{2}\leq 2\mE\int^{t}_{t-\delta}(\delta|f(x(r), \mu_{r})
 +u(x(\sigma_{r}), \mu_{\sigma_{r}})|^{2}+|g(x(r), \mu_{r})|^{2})\dif r.
\end{align*}
Choosing $\theta\geq \frac{(\frac{1}{\lambda_{1}}+\frac{1}{\lambda_{2}})L_{3}}{1-8L_{3}\delta^{2}},$  this together with \eqref{16} yields that
\begin{align}\label{18}
 \mE \mathscr{L}{V}(x_{t}, \mu_{t})\leq -\lambda_{4}\mE|x(t)|^{2}+\gamma_{2}.
\end{align}
From \eqref{8}, we get
\begin{align*}
 0\leq \mE \widetilde{V}(x_{t}, \mu_{t})\leq L_{4}- \lambda_{4}\mE\int^{t}_{0}|x(s)|^{2}\dif s+\gamma_{2} t.
\end{align*}
where $L_{4}=V(x_{0},\mu_{0})+\theta\int^{0}_{-\delta}\int^{0}_{r}[\delta|f(x(u), \mu_{u})
 +u(x(\sigma_{u}), \mu_{\sigma_{u}})|^{2}+|g(x(u), \mu_{u})|^{2}]\dif u\dif r.$

This   leads to
$$\limsup_{t\rightarrow \infty}\frac{1}{t}\mE\int^{t}_{0}|x(s)|^{2}\dif s\leq \gamma_{2}. $$
If $\gamma_{2}=0, $ the second assertion follows.
\end{proof}

\bl\label{13.2y}
Assume that  $(\mathrm{H1})$ and  $(\mathrm{H2})$ hold. Let $\delta > 0$ be sufficiently small such that  $H(\delta, p):=(3^{p-1}\delta^{2} 2^{L_{1}+\frac{p}{2}}+3^{p-1}c_{p}\delta^{2} 2^{L_{2}+\frac{p}{2}}+3^{p-1}\delta^{2} 2^{L_{2}+\frac{p}{2}})e^{3^{p-1}\delta 2^{L_{1}+\frac{p}{2}}+3^{p-1}\delta c_{p} 2^{L_{2}+\frac{p}{2}}}$ $<\frac{1}{2^{p}},$  where $c_{p}$ is the constant in BDG's inequality.  Then the solution $x(t)$ of Eq.\eqref{2} satisfies the following inequality for $p\geq 2$:
\begin{align}\label{20}
 \mE|x(t)-x(\sigma_{t})|^{p}\leq \frac{2^{p-1}H(\delta,p)}{1-2^{p-1}H(\delta,p)}\mE|x(t)|^{p}.
\end{align}
\el
\begin{proof}
Fix any integer $l\geq 0.$ For $t\in [l\delta, (l+1)\delta),$  we have $\sigma_{t}=\lfloor \frac{t}{\delta} \rfloor \delta
=l\delta.$
From \eqref{2}, we obtain
\begin{align*}
 |x(t)-x(\sigma_{t})|^{p}&= |x(t)-x(l\delta)|^{p}\no\\
 &=\bigg|\bigg.\int^{t}_{l\delta}(f(x(s),\mu_{s})+u(x(\sigma_s),\mu_{\sigma_s})) \dif s + \int^{t}_{l\delta}g(x(s),\mu_{s}) \dif B(s)  \bigg|\bigg.^{p}\no\\
 &\leq 3^{p-1}\bigg|\bigg.\int^{t}_{l\delta}f(x(s),\mu_{s}) \dif s\bigg|\bigg.^{p} + 3^{p-1}\bigg|\bigg.\int^{t}_{l\delta}u(x(\sigma_s),\mu_{\sigma_s}) \dif s\bigg|\bigg.^{p}\no\\
 &+ 3^{p-1}\bigg|\bigg.\int^{t}_{l\delta}g(x(s),\mu_{s}) \dif B(s)  \bigg|\bigg.^{p}.
\end{align*}
This, together with the  Lipschitz condition $(\mathrm{H1})$, implies
\begin{align*}
 \mE&|x(t)-x(\sigma_{t})|^{p}\no\\
&\leq 3^{p-1}\mE\bigg|\bigg.\int^{t}_{l\delta}f(x(s),\mu_{s}) \dif s\bigg|\bigg.^{p} + 3^{p-1}\mE\bigg|\bigg.\int^{t}_{l\delta}u(x(\sigma_s),\mu_{\sigma_s}) \dif s\bigg|\bigg.^{p}\no\\
 &+ 3^{p-1}\mE\bigg|\bigg.\int^{t}_{l\delta}g(x(s),\mu_{s}) \dif B(s)  \bigg|\bigg.^{p}\no\\
 &\leq (3^{p-1}\delta 2^{L_{1}+\frac{p}{2}}+3^{p-1}c_{p}\delta 2^{L_{2}+\frac{p}{2}})\mE\int^{t}_{l\delta}|x(s)-x(l\delta)|^{2}\dif s\no\\
 &+(3^{p-1}\delta^{2} 2^{L_{1}+\frac{p}{2}}+3^{p-1}\delta^{2} 2^{L_{2}+\frac{p}{2}}+3^{p-1}c_{p}\delta^{2} 2^{L_{2}+\frac{p}{2}})\mE|x(l\delta)|^{p}.
\end{align*}
It follows from Gronwall's inequality  that
\begin{align*}
 \mE|x(t)-x(\sigma_{t})|^{p}
\leq  H(\delta,p)\mE|x(\sigma_{t})|^{p}.
\end{align*}
Hence, the required assertion follows from $2^{p-1}H(\delta,p)<1$ and
\begin{align*}
 \mE|x(t)-x(\sigma_{t})|^{p}
\leq  2^{p-1}H(\delta,p)(\mE|x(t)|^{p}+\mE|x(t)-x(\sigma_{t})|^{p}).
\end{align*}
\end{proof}
The following theorem states the asymptotic stability in mean square of the solution of Eq. \eqref{2}.
\bt\label{T3.3}
Assume $(\mathrm{H1})$ and    $(\mathrm{H2})$  hold with $\gamma_{2}=0$. If $\delta>0$ is sufficiently small  such that
$H(\delta):=(12\delta^{2} L_{1}+6\delta^{2}L_{3}+12\delta L_{2})e^{(12\delta L_{1}+12L_{2})\delta}<\frac{1}{2},$  then  the solution of controlled system \eqref{2} is stable in mean square, i.e.,
$$\lim_{t\rightarrow \infty}\mE|x(t)|^{2}=0.$$
\et
\begin{proof}
For $0\leq t_{1}< t_{2}<\infty,$  we have
\begin{align*}
 &\mE|x(t_{2})-x(t_{1})|^{2}\no\\
&\leq 2|t_{2}-t_{1}|\mE\int^{t_{2}}_{t_{1}}|f(x(s),\mu_{s})+u(x(\sigma_s),\mu_{\sigma_s})|^{2}\dif s+2\mE\int^{t_{2}}_{t_{1}}|g(x(s),\mu_{s})|^{2}\dif s\no\\
&\leq 4|t_{2}-t_{1}|\mE\int^{t_{2}}_{t_{1}}(|f(x(s),\mu_{s})|^{2}+|u(x(\sigma_s),\mu_{\sigma_s})|^{2})\dif s+2\mE\int^{t_{2}}_{t_{1}}|g(x(s),\mu_{s})|^{2}\dif s\no\\
&\leq 16(L_{1}+2L_{3})|t_{2}-t_{1}|\int^{t_{2}}_{t_{1}}\mE|x(s)|^{2}\dif s+8L_{2}\int^{t_{2}}_{t_{1}}\mE|x(s)|^{2}\dif s\no\\
&+32L_{3}|t_{2}-t_{1}|\int^{t_{2}}_{t_{1}}\mE|x(s)- x(\sigma_{s})|^{2}\dif s\no\\
&\leq 16\bigg(\bigg. L_{1}+2L_{3}+2L_{3}\frac{H(\delta)}{1-2H(\delta)}\bigg)\bigg.|t_{2}-t_{1}|\int^{t_{2}}_{t_{1}}\mE|x(s)|^{2}\dif s+8L_{2}\int^{t_{2}}_{t_{1}}\mE|x(s)|^{2}\dif s.
\end{align*}
By \eqref{91}, one can see  that
 $\int^{t}_{0}\mE|x(s)|^{2}\dif s$ is uniformly continuous in $t$ on $\mR^{+}.$
Therefore, $\mE|x(t)|^{2}$ is uniformly continuous in $t.$   This together with \eqref{91}  implies the assertion.
\end{proof}

Next, we will present the  exponential stability in mean square.
\bt
 Assume that $(\mathrm{H1})$   and  $(\mathrm{H2})$  hold with $\gamma_{2}=0$. Let $\delta$ be sufficiently small such that  $\lambda_{4}=\gamma_{1}c_{1}-8L_{1}\theta\delta^{2}-8L_{3}\theta\delta^{2}-2L_{2}\theta\delta > 0$. If there exist two positive constants $c_{1}$ and  $c_{2}$ such that
\begin{align}\label{29}
c_{1}\int_{\mR^{d}}|x|^{2}\mu(\dif x)\leq \int_{\mR^{d}}V(x,\mu)\mu(\dif x)\leq c_{2}\int_{\mR^{d}}|x|^{2}\mu(\dif x),
\end{align}
then the  solution $x(t)$ of Eq. \eqref{2} is exponentially stable in mean square, i.e.,
\begin{align*}
 \limsup_{t\rightarrow \infty}\frac{1}{t}\log(\mE|x(t)|^{2})\leq -\alpha,
\end{align*}
where   $\alpha>0$ is a constant satisfying $\alpha K\tau e^{\alpha\delta}+\alpha c_{2}-\lambda_{4}<0 $ with $K=\frac{8\theta\delta^{2}L_{3}H(\delta)}{1-2H(\delta)}+4\theta \delta^{2}L_{1}+8\theta \delta^{2}L_{3}+2\theta\delta L_{2}.$
\et
\begin{proof}
Let $\widetilde{V}(x_{t}, \mu_{t})$ be defined by \eqref{7}. Applying  It\^{o}'s formula to $e^{\alpha t}\widetilde{V}(x_{t}, \mu_{t})$ and using (\ref{8}), we get
\begin{align*}
\mE e^{\alpha t}\widetilde{V}(x_{t}, \mu_{t})=\mE\widetilde{V}(x_{0}, \delta_{x_{0}}) +\mE\int^{t}_{0}e^{\alpha s}[\alpha \widetilde{V}(x_{s}, \mu_{s})+\mathscr{L}{V}(x_{s},\mu_{s})]\dif s.
\end{align*}
By \eqref{18}, we obtain
\begin{align}\label{32}
\mE e^{\alpha t}\widetilde{V}(x_{t}, \mu_{t})\leq \mE\widetilde{V}(x_{0}, \delta_{x_{0}}) +\int^{t}_{0}e^{\alpha s}[\alpha \mE\widetilde{V}(x_{s}, \mu_{s})-\lambda_{4}\mE|x(s)|^{2}]\dif s.
\end{align}
Due to  \eqref{7} and \eqref{29}, one can see that
\begin{align*}
\mE \widetilde{V}(x_{t}, \mu_{t})\leq c_{2}\mE|x(t)|^{2}+\mE(\Gamma(x_{t}, \mu_{t})),
\end{align*}
where $\Gamma(x_{t}, \mu_{t})=\theta\int^{t}_{t-\delta}\int^{t}_{r}[\delta|f(x(u), \mu_{u})
 +u(x(\sigma_{u}), \mu_{\sigma_{u}})|^{2}+|g(x(u), \mu_{u})|^{2}]\dif u\dif r.$
It follows  from $(\mathrm{H1})$ that
\begin{align*}
\mE \Gamma(x_{t}, \mu_{t})&\leq 4\theta\delta^{2}L_{3}\mE\int^{t}_{t-\delta}|x(s)-x(\sigma_{s})|^{2}\dif s \no\\
&+4\theta\delta^{2}L_{3}\mE\int^{t}_{t-\delta}W_{2}^{2}(\mu_s,\mu_{\sigma_{s}})\dif s \no\\
&+(2\theta \delta^{2}L_{1}+4\theta \delta^{2}L_{3}+\theta\delta L_{2})\mE\int^{t}_{t-\delta}|x(s)|^{2}\dif s\no\\
&+(2\theta \delta^{2}L_{1}+4\theta\delta^{2}L_{3}+\theta\delta L_{2})\mE\int^{t}_{t-\delta}W_{2}^{2}(\mu_{s},\delta_{0})\dif s\no\\
&\leq 8\theta\delta^{2}L_{3}\mE\int^{t}_{t-\delta}|x(s)-x(\sigma_{s})|^{2}\dif s \no\\
&+(4\theta \delta^{2}L_{1}+8\theta \delta^{2}L_{3}+2\theta\delta L_{2})\mE\int^{t}_{t-\delta}|x(s)|^{2}\dif s.
\end{align*}
By Lemma 3.3, we have
\begin{align*}
\mE \Gamma(x_{t}, \mu_{t})
&\leq \frac{8\theta\delta^{2}L_{3}H(\delta,2)}{1-2H(\delta,2)}\mE\int^{t}_{t-\delta}|x(s)|^{2}\dif s \no\\
&+(4\theta \delta^{2}L_{1}+8\theta \delta^{2}L_{3}+2\theta\delta L_{2})\mE\int^{t}_{t-\delta}|x(s)|^{2}\dif s\no\\
&\leq K\mE\int^{t}_{t-\delta}|x(s)|^{2}\dif s.
\end{align*}
From  \eqref{29} and \eqref{32}, we derive that
\begin{align*}
c_{1}e^{\alpha t}\mE |x(t)|^{2}&\leq \widetilde{V}(x_{0}, \delta_{x_{0}}) +  (\alpha K\delta e^{\alpha\delta}+\alpha c_{2}-\lambda_{4})\int^{t}_{0}\mE|x(s)|^{2}\dif s.
\end{align*}
This together with $\alpha K\delta e^{\alpha\delta}+\alpha c_{2}-\lambda_{4}<0$ yields
\begin{align*}
c_{1}e^{\alpha t}\mE |x(t)|^{2}\leq \widetilde{V}(x_{0}, \delta_{x_{0}}), \, \, t>\delta.
\end{align*}
 The required assertion follows.
\end{proof}

\section{Interacting particle systems }

Assume that $\{B^{1}(t)\}, \{B^{2}(t) \},  \cdots,$  are independent $m-$dimensional Brownian motions. We now consider the following equations, for $i=1,2,\cdots$

\begin{equation}\label{2++}
\begin{cases}
 &\dif x^{i}(t)= (f( x^{i}(t), \mu^{x_{i}}_{t})+u( x^{i}(\sigma_{t}), \mu^{x_{i}}_{\sigma_{t}}))\dif t +g(x^{i}(t), \mu^{x_{i}}_{t})\dif  B^{i}(t), t\in [0, \infty),\\
 & x(0)=x_{0},
 \end{cases}
 \end{equation}
where $\mu^{x_{i}}_{t}$ represents the law of $x_{i}(t).$  Let $\{x^{i}(t), i=1,2,\cdots\}$ be the unique solution   of the above equations.
We now  write the corresponding interacting particle systems as follows:

 \begin{equation}\label{2+}
\begin{cases}
 &\dif x^{i,N}(t)= (f( x^{i,N}(t), \mu^{x,N}_{t})+u( x^{i,N}(\sigma_{t}), \mu^{x,N}_{\sigma_{t}}))\dif t +g(  x^{i,N}(t), \mu^{x,N}_{t})\dif  B^{i}(t), t\in [0, \infty),\\
 & x(0)=x_{0},
 \end{cases}
 \end{equation}
where $\mu^{x,N}_{t}(\cdot):= \frac{1}{N}\sum^{N}_{i=1}\delta_{x^{i,N}(t)}(\cdot).$
Obviously, in real world, the distribution of  $x(t)$ is  difficult to be observed. However,    the corresponding one of the particle system can be observed.  We will prove that the exponential stability of system \eqref {2} is equivalent to   the  the exponential stability of  corresponding particle system \eqref {2+}.

 First of all,  we make the following assumption:
\begin{itemize}
\item[(H3)] Let   $\xi$ and $\eta$ be  two random variables  whose distributions  are $\mu$ and $\nu$,  respectively, and  the joint distribution of  $(\xi, \eta)$ be $F_{\xi, \eta}(z, \bar z).$   Assume that there exists a Lyapunov  function $U \in  \sC(\mR^{d}\times\sM_{\lambda^{2}}(\mR^{d}))$ such that
\begin{align*}
\mE&\int_{\mR^{d}}\int_{\mR^{d}}\mathbb{L}{U}(x, \mu,y, \nu)F_{\xi, \eta}(x, y)\no\\
&\leq -\alpha_{1} \int_{\mR^{d}}|x|^{p}\mu(\dif x) +\alpha_{2}W^{p}_{2}(\mu, \nu) +\alpha_{3}\int_{\mR^{d}}\int_{\mR^{d}}|x-y|^{p}F_{\xi, \eta}(x, y) +\beta
\end{align*}
for $p\geq2,$ where $\alpha_{1},\alpha_{2},\alpha_{3},\beta$ are four constants satisfying $\alpha_{1}>0,\alpha_{2}\geq 0, \alpha_{3}\geq 0, \beta\geq 0. $
\end{itemize}

For the future use, we cite  \cite[Theorem 5.8, pp.362]{CD}  as the following lemma.
\bl\label{L4.1Y}
Assume that  $\{x_{n}\}_{n\geq 1}$  is a sequence of independent identically  distributed (i.i.d.  for short) random variables in $\mR^{d}$ with common distribution $\mu \in \sP(\mR^{d}).$ For any $N\in \mN,$  we define the empirical measure $\mu^{N}= \frac{1}{N}\sum^{N}_{i=1}\delta_{x_{i}}.$  If $\mu\in \sP_{q}(\mR^{d})$  with $q> 4,$  then there exists a constant $C=C(d,q, W_{q}(\mu))$ such that for any $N\geq 2,$
\begin{align*}
\mE[W^{2}_{2}(\mu^{N}, \mu)]\leq\tau(N)= C
 \begin{cases}
N^{-\frac{1}{2}}, &1\leq d< 4,\\
N^{-\frac{1}{2}}\ln(N), &d= 4,\\
N^{-\frac{2}{ d}}, &4< d.
\end{cases}
\end{align*}
\el

The constant $C$ in the lemma above depends on the $q$th moment  $W_q(\mu)$ of i.i.d. random variables. In order to apply this result to  the solution $x(t)$ of Eq. \eqref{2}, we give the following moment estimate of $x(t)$.

\bl\label{L4.3Y}
Assume $\mathrm{(H1)},\mathrm{(H2)}$ and $\mathrm{(H3)}.$ 
Then it holds  that
\begin{align*}
\sup_{t\geq 0}\mE|x(t)|^{p}\leq C_{x_0, p},\quad  p \ge 2.
\end{align*}
where   $C_{x_0, p}$ only depends on $x_{0}, p.$
\begin{proof}
Let $\alpha, \delta$ be  two positive constants sufficiently small such that
$$H(\delta, p)< \frac{1}{2}\,\,\mbox{and}\,\,\alpha-\alpha_{1} +\alpha_{2}\frac{H(\delta,p)}{1-2H(\delta,p)}+\alpha_{3}\frac{H(\delta,p)}{1-2H(\delta,p)}< 0. $$
From It\^{o}'s formula, $(\mathrm{H3})$ and \eqref{20}, we have
\begin{align*}
&\mE[e^{\alpha t}U(x(t), \mu_{t})]= U(x_0, \mu_{0})
 + \mE\int^{t}_{0}e^{\alpha s}(\alpha +\mathbb{L}{U}(x(s), \mu_{s},x(\sigma_{s}), \mu_{\sigma_{s}}))\dif s\no\\
& \leq   U(x_0, \mu_{0})
 + \mE\int^{t}_{0}e^{\alpha s}(\alpha-\alpha_{1} |x(s)|^{p} +\alpha_{2}W^{p}_{p}(\mu_{s}, \delta_{0}) +\alpha_{3}|x(s)-x(\sigma_{s})|^{p} +\beta)\dif s\no\\
&  \leq U(x_0, \mu_{0})
 + \mE\int^{t}_{0}e^{\alpha s}\bigg[\bigg.\bigg(
 \bigg.\alpha-\alpha_{1} +\alpha_{2}\frac{H(\delta,p)}{1-2H(\delta,p)}+\alpha_{3}\frac{H(\delta,p)}{1-2H(\delta,p)}\bigg)\bigg.|x(s)|^{p}+\beta)\bigg]\bigg.\dif s.
 \end{align*}
Due to the assumption in the theorem,   we have
\begin{align*}
 \mE[|x(t)|^{p}]\leq e^{-\alpha t}U(x_0, \mu_{0})+ \frac{\beta}{\alpha}(1-e^{-\alpha t}).
\end{align*}
 We obtain the required results from the above inequality.
\end{proof}

\el
From the above lemma, one can see that $\sup_{t\geq 0} W_{q}(\mu_{t})\leq C_{q}.$  Thus we have the following theorem.
\bt
Assume $\mathrm{(H1)}-\mathrm{(H3)}$, and $p> 4.$  Then, we have
\begin{align*}
\sup_{0\leq i\leq N}\sup_{t\geq 0}\mE[|x^{i}(t)- x^{i,N}(t)|^{2}]
 \leq C_{N}:=C\begin{cases}
N^{-\frac{1}{2}}, &1\leq d< 4,\\
N^{-\frac{1}{2}}\ln(N), &d= 4,\\
N^{-\frac{2}{ d}}, &4< d,
\end{cases}
\end{align*}
where $C$ only depends on $d,q, \sup_{t\geq 0}\mE|x(t)|^{p}.$
\et
\begin{proof}
By It\^{o}'s formula and Assumption (H1),  we derive
\begin{align}\label{40}
\mE|x^{i}(t)-x^{i,N}(t)|^{2}&\leq \mE\bigg\{\bigg. \int^{t}_{0}\bigg[\bigg.2\langle x^{i}(s)-x^{i,N}(s), f(x^{i}(s),\mu^{x^{i}}_{s})-f(x^{i,N}(s), \mu_{s}^{x,N}) \rangle \no\\
&+ |g(x^{i}(s), \mu^{x^{i}}_{s})-g(x^{i,N}(s), \mu_{s}^{x,N})|^{2}  \bigg]\bigg. \dif s     \bigg\}\bigg.\no\\
&\leq L\int^{t}_{0}\mE|x^{i}(s)-x^{i,N}(s)|^{2}\dif s +L\int^{t}_{0}\mE[W^{2}_{2}(\mu^{x^{i}}_{s},\mu_{s}^{x,N})]\dif s,
\end{align}
where $L$ is a constant being independent of $t,$ whose value  may vary from one place to another.  We construct another empirical measure which comes from \eqref{2++} as follows:
$$\mu^{N}_{s}(\dif x)= \frac{1}{N}\sum^{N}_{j=1}\delta_{x^{j}(s)}(\dif x).$$
Note that
\begin{align*}
W^{2}_{2}(\mu^{x^{i}}_{s},\mu_{s}^{x,N})&\leq 2(W^{2}_{2}(\mu^{x^{i}}_{s},\mu^{N}_{s})+W^{2}_{2}(\mu^{N}_{s},\mu_{s}^{x,N}))\no\\
& = 2W^{2}_{2}(\mu^{x^{i}}_{s},\mu^{N}_{s})+\frac{2}{N}\sum^{N}_{j=1}\mE|x^{j}(s)-x^{j,N}(s)|^{2}.
\end{align*}
and
\begin{align*}
\frac{1}{N}\sum^{N}_{j=1}\mE|x^{j}(s)-x^{j,N}(s)|^{2}= \mE|x^{j}(s)-x^{j,N}(s)|^{2}.
\end{align*}
Thus, we have
\begin{align*}
\int^{t}_{0}\mE[W^{2}_{2}(\mu^{x^{i}}_{s},\mu_{s}^{x,N})]\dif s
\leq 2\int^{t}_{0}\mE[W^{2}_{2}(\mu^{x^{i}}_{s},\mu^{N}_{s})]\dif s +2\int^{t}_{0}\mE|x^{j}(s)-x^{j,N}(s)|^{2}\dif s.
\end{align*}
This together with \eqref{40} implies
\begin{align}\label{i2y}
\mE|x^{i}(t)-x^{i,N}(t)|^{2}
&\leq L\int^{t}_{0}\mE|x^{i}(s)-x^{i,N}(s)|^{2}\dif s +L\int^{t}_{0}\mE[W^{2}_{2}(\mu^{x^{i}}_{s},\mu^{N}_{s})]\dif s.
\end{align}
By Lemma \ref{L4.1Y}, one has
$$
\mE[W^{2}_{2}(\mu^{x^{i}}_{s},\mu^{N}_{s})]\le C(d, q, W_q(\mu^i))\begin{cases}
N^{-\frac{1}{2}}, &1\leq d< 4,\\
N^{-\frac{1}{2}}\ln(N), &d= 4,\\
N^{-\frac{2}{ d}}, &4< d,
\end{cases}
$$
where $\mu^i$ is the distribution of $x^{i}(s)$.  Additionally, Lemma \ref{L4.3Y} implies that $\sup_{s \ge 0}C(d, q, W_q(\mu^i))<\infty$.

Therefore, the proof is complete by Gronwall's inequality and \eqref{i2y}.
\end{proof}
\bt
 The solution  of Eq. \eqref{2}  is  exponentially stable in mean square, i.e.,  there exists a positive constant $\ell_{1}$ such that
 \begin{align}
 \limsup_{t\rightarrow \infty}\frac{1}{t}\log(\mE|x(t)|^{2})\leq -\ell_{1},
\end{align}
if and only if the solution of \eqref{2+} is  exponentially stable in mean square, i.e.,  there exists a positive constant $\ell_{2}$ and for any $i$ such that
 \begin{align}
 \limsup_{t\rightarrow \infty}\lim_{N\rightarrow \infty}\frac{1}{t}\log(\mE|x^{i,N}(t)|^{2})\leq -\ell_{2}.
\end{align}
\et

\begin{proof}  One may complete the proof by the following two inequalities:
\begin{align*}
 \mE|x^{i,N}(t)|^{2}\leq 2\mE|x^{i}(t)-x^{i,N}(t)|^{2} +2\mE|x^{i}(t)|^{2}
 \leq 2C_{N}+2\mE|x(t)|^{2},\\
 \mE|x^{i}(t)|^{2}\leq 2\mE|x^{i}(t)-x^{i,N}(t)|^{2} +2\mE|x^{i, N}(t)|^{2}
 \leq 2C_{N}+2\mE|x^{i, N}(t)|^{2}.
\end{align*}
\end{proof}

We now give an example to illustrate the theory.
\begin{exa}
 Consider the following equation:
 {\rm \begin{equation}\label{u}
\begin{cases}
&\dif y(t)=\bigg(\bigg.2y(t)+\int_{\mR}z\mu_{t}(\dif z)\bigg)\bigg.\dif t + y(t)\dif B(t),\\
&y(0)=y_{0},
\end{cases}
\end{equation}
where $y_{0}$ is a positive constant. Setting  $V(x, \mu)=|x|^{2}+\int_{\mR}|z|^{2}\mu(\dif z), $  by the fact of $\partial_{\mu}\bigg(\bigg.\int_{\mR}|z|^{2}\mu(\dif z)\bigg)\bigg.(y)=2y,$ we have
\begin{align*}
LV(x,\mu)&=(2x+\int_{\mR}z\mu(\dif z))2x+ |x|^{2}+\int_{\mR}\bigg(\bigg.2y+\int_{\mR}z\mu(\dif z)\bigg)\bigg.2y\mu(\dif y)+\int_{\mR}y^{2}\mu(\dif y)\no\\
&\geq 4|x|^{2}+5\int_{\mR}|z|^{2}\mu(\dif z)+\bigg|\bigg.\int_{\mR}z\mu(\dif z) \bigg|\bigg.^{2},
\end{align*}
where $k_{1},  k_{2}$ are two undetermined constants which will be given in the following.
Thus, from It\^{o}'s formula,  we can  know that
the solution of Eq.\eqref{u} is unstable in the sense of mean square expectation.

We now consider the following equation with discrete time feedback control:
\begin{align}\label{s}
\dif x(t)=[2x(t)+\int_{\mR}z\mu_{t}(\dif z) -k_{1}x(\sigma_{t})-k_{2}\int_{\mR}z\mu_{\sigma_{t}}(\dif z)]\dif t + x(t)\dif B(t),
\end{align}
 where $k_1$ and $k_2$ are constants. Computing the operator of  Eq.\eqref{s} acting on $V(x, \mu)$, one can see that
\begin{align*}
LV(x,\mu)&=\bigg(\bigg.2x+\int_{\mR}z\mu(\dif z)-k_{1}x-k_{2}\int_{\mR}z\mu(\dif z)\bigg)\bigg.2x+ |x|^{2}\no\\
&+\int_{\mR}\bigg(\bigg.2y+\int_{\mR}z\mu(\dif z)-k_{1}y-k_{2}\int_{\mR}z\mu(\dif z)\bigg)\bigg.2y\mu(\dif y)+\int_{\mR}y^{2}\mu(\dif y)\no\\
&\leq (6-2k_{1}+k_{2})|x|^{2}+(5-2k_{1})\int_{\mR}|z|^{2}\mu(\dif z)+(3-k_{2})\bigg(\bigg.\int_{\mR}z \mu(\dif z)   \bigg)\bigg.^{2}.
\end{align*}
and
\begin{align*}
\lambda_{1}\int_{\mR}|V_{x}(x, \mu)|^{2}\mu(\dif x)=4\lambda_{1}\int_{\mR}|x|^{2}\mu(\dif x), \lambda_{2}\int_{\mR}|\partial_{\mu}V_{x}(x, \mu)|^{2}\mu(\dif x)\leq 4\lambda_{2}\int_{\mR}|x|^{2}\mu(\dif x).
\end{align*}
 Choosing $\lambda_{1}=\frac{1}{2}, \lambda_{2}=\frac{1}{2}, k_{1}=8, k_{2}=3,$  we have
\begin{align}\label{45}
\int_{\mR}LV&(x,\mu)\mu(\dif x)+\lambda_{1}\int_{\mR}|V_{x}(x, \mu)|^{2}\mu(\dif x)+\lambda_{2}\int_{\mR}|\partial_{\mu}V_{x}(x, \mu)|^{2}\mu(\dif x)
\no\\
& \leq  -7\int_{\mR}x^{2}\mu(\dif x).
\end{align}
Obviously,  $(\mathrm{H1})$ holds,  and $(\mathrm{H2})$ holds with $\gamma_{1}= -5, \gamma_{2}=0.$   Moreover, $\int_{\mR}|x|^{2}\mu(\dif x)\leq\int_{\mR}V(x, \mu)\mu(\dif x)\leq 2\int_{\mR}|x|^{2}\mu(\dif x). $  This means that  the conditions of Theorem 3.5 hold. Therefore, we conclude that  the solution of Eq.\eqref{s} is   exponentially stable in mean square.
Set $U(x, \mu)=|x|^{6}+\int_{\mR}|z|^{6}\mu(\dif z).$  Furthermore,   from   \eqref{12}  and Lemma \ref{13.2y}, we know that
 \begin{align}\label{e}
&\mE[\mathbb{L}U(x(t), \mu_t, x(\sigma_t), \mu_{\sigma_t})]\no\\
&=\mE[LU(x(t), \mu_{t})-\partial_{x}U(x(t), \mu_{t})(u(x(t), \mu_{t})-u(x(\sigma_{t}), \mu_{\sigma_{t}}))] \no\\
 &- \int_{\mR^{d}}\int_{\mR^{d}}\partial_{\mu}V(x(t), \mu_{t})(y)(u(y,\mu_{t})-u(\bar{y},\mu_{\sigma_t}))F_{x(t), x(\sigma(t)) }(\dif y, \dif \bar{y})\no\\
 &\leq (52-12k_{1}+12k_{2}+\frac{5}{3})\mE[|x(t)|^{6}] +6^{5}2L^{3}_{3}\mE[|x(t)-x(\sigma_{t})|^{6}]\no\\
 &\leq (52-12k_{1}+12k_{2}+\frac{5}{3})\mE[|x(t)|^{6}] + 6^{5}2L^{3}_{3}\frac{32H(\delta,6)}{1-32H(\delta,6)}\mE[|x(t)|^{6}]  \no\\
 &\leq -\frac{13}{3}\mE[|x(t)|^{6}]+6^{5}2L^{3}_{3}\frac{32H(\delta,6)}{1-32H(\delta,6)}\mE[|x(t)|^{6} ].
 \end{align}
  Letting $\delta$  be small enough such that $6^{5}2L^{3}_{3}\frac{32H(\delta,6)}{1-32H(\delta,6)}< 1,$  we infer that the conditions of Lemma \ref{L4.3Y} hold.  Thus, the corresponding interacting particle system is exponentially stable in mean square.

 }
\end{exa}

\section*{Funding}
This research is supported by the National Natural Science Foundation of China (Grant no. 61876192, 11626236), the Fundamental Research Funds for the Central Universities of South-Central University for Nationalities (Grant nos. CZY15017, KTZ20051, CZT20020).

\end{document}